\documentclass   [12pt]{article}          

\begin{document}



\catcode`@=11 \catcode`!=11

\expandafter\ifx\csname fiverm\endcsname\relax
  \let\fiverm\fivrm
\fi
  
\let\!latexendpicture=\endpicture 
\let\!latexframe=\frame
\let\!latexlinethickness=\linethickness
\let\!latexmultiput=\multiput
\let\!latexput=\put
 
\def\@picture(#1,#2)(#3,#4){%
  \@picht #2\unitlength
  \setbox\@picbox\hbox to #1\unitlength\bgroup 
  \let\endpicture=\!latexendpicture
  \let\frame=\!latexframe
  \let\linethickness=\!latexlinethickness
  \let\multiput=\!latexmultiput
  \let\put=\!latexput
  \hskip -#3\unitlength \lower #4\unitlength \hbox\bgroup}

\catcode`@=12 \catcode`!=12
\input{pictex.tex}

 
\catcode`@=11 \catcode`!=11
  
\let\!pictexendpicture=\endpicture 
\let\!pictexframe=\frame
\let\!pictexlinethickness=\linethickness
\let\!pictexmultiput=\multiput
\let\!pictexput=\put

\def\beginpicture{%
  \setbox\!picbox=\hbox\bgroup%
  \let\endpicture=\!pictexendpicture
  \let\frame=\!pictexframe
  \let\linethickness=\!pictexlinethickness
  \let\multiput=\!pictexmultiput
  \let\put=\!pictexput
  \let\input=\@@input   
  \!xleft=\maxdimen  
  \!xright=-\maxdimen
  \!ybot=\maxdimen
  \!ytop=-\maxdimen}

\let\frame=\!latexframe

\let\pictexframe=\!pictexframe

\let\linethickness=\!latexlinethickness
\let\pictexlinethickness=\!pictexlinethickness

\let\\=\@normalcr
\catcode`@=12 \catcode`!=12


\newcommand{\bbeta}{\mbox{\boldmath $\beta$}} 
\newcommand{\BO}{\mbox{\rm O}}                
\newcommand{\m}{\mbox{\boldmath $m$}}         
\newcommand{\M}{\mbox{\boldmath $M$}}         
\newcommand{\mi}{{$-$}}                       
\newcommand{\n}{\mbox{\boldmath $n$}}         
\newcommand{\p}{\phantom{5}}                  
\newcommand{\ps}{\phantom{$^*$}}              
\newcommand{\pc}{{\scriptsize$\bullet$}}      
\newcommand{\q}{\phantom{55}}                 
\newcommand{\qs}{\hspace*{.1in}}              
\newcommand{\rme}{{\rm e}}                    
\newcommand{\U}{\mbox{\boldmath $u$}}         
\newcommand{\X}{\mbox{\boldmath $X$}}         
\newcommand{\x}{\mbox{\boldmath $x$}}         
\newcommand{\z}{\mbox{\boldmath $z$}}         

\hyphenation{general-iza-tions}
\hyphenation{distri-bution}
\hyphenation{Zel-ter-man}


\thispagestyle{empty}

\begin{center}
\vspace*{\fill}

{\Large\bf               The Maximum Negative            }\\[1ex]
{\Large\bf            Hypergeometric Distribution        }\\

\vspace*{1.5in}
\begin{tabular}{l}
{\bf Daniel Zelterman}                \\[.2in]

Division of Biostatistics             \\
School of Medicine                    \\
Yale University                       \\
New Haven, CT  06520                  \\[.4in]

\today
\end{tabular}
\end{center}

\vspace*{\fill}

\noindent
Email address for author:
{\tt daniel.zelterman@yale.edu} . 
This research was supported by grant P30-CA16359 awarded by the U.S.~National Institutes of Health to the Yale Comprehensive Cancer Center.


\newpage
\vspace*{.25in}
\thispagestyle{empty}

\section*                      {\bf   Abstract}

An urn contains a known number of balls of two different colors.
We describe the random variable counting the smallest number of draws needed in order to observe at least $\,c\,$ of both colors when sampling without replacement for a prespecified value of $\,c=1,2,\ldots\,$.
This distribution is the finite sample analogy to the maximum negative binomial distribution described by Zhang, Burtness, and Zelterman~(2000).
We describe the modes, approximating distributions, and estimation of the contents of the urn.

\vspace*{.25in}

\noindent {\bf Keywords:} discrete distributions;
negative binomial distribution;
riff-shuffle distribution;
hypergeometric distribution; 
negative hypergeometric distribution;
maximum negative binomial distribution

\vspace*{.25in}


\thispagestyle{empty}
\setcounter{page}{1}
\section                          {Introduction}

{\centering\it

\noindent And the LORD said unto Noah, Come thou and all thy house into the ark; 

for thee have I seen righteous before me in this generation.

\quad Of every clean beast thou shalt take to thee by sevens, the male and his female: 

\quad\quad\quad\quad\quad\quad and of beasts that are not clean by two, the male and his female.}

Genesis 7:1--2.  King James translation

This charge to Noah required seven pairs of clean animals.
How many animals did Noah plan on catching in order to be reasonably sure of achieving $\,c=7\,$ male and female pairs? 
He didn't want to handle more dangerous, wild creatures than necessary.
In the case of a rare or endangered species, the finite population size $\,N\,$ could be small.
A ``clean'' animal meant it was suitable for consumption or sacrifice.

In a sequence of independent and identically distributed Bernoulli~$(p)\,$ random variables, the {\em negative binomial distribution\/} describes the behavior of the number of failures $\,Y\,$ observed before observing $\,c\,$ successes, for integer valued parameter $\,c\geq 1$.
This well-known distribution has probability mass function
\begin{equation}                                                  
   \Pr[\,Y=y\,] = {{c+y-1}\choose{c-1}}p^c (1-p)^y              \label{NB}
\end{equation}
defined for $\,y=0,1,\ldots\,$.

The negative binomial distribution~(\ref{NB}) is discussed in detail by Johnson, Kotz, and Kemp (1992, Ch.~5).
In this introductory section we will describe several sampling schemes closely related to the negative binomial.
Table~1 may be useful in illustrating the various relations between these distributions.

The {\em maximum negative binomial distribution\/} is the distribution of the smallest number of trials needed in order to observe at least $\,c\,$ successes and $\,c\,$ failures for integer valued parameter $\,c\geq 1$.
This distribution is motivated by the design of a medical trial in which we want to draw inference on the Bernoulli parameter $\,p\,$ in an infitely large population.
If the prevalance $\,p\,$ of a binary valued genetic trait in cancer patients is very close to either zero or one then there is little to be gained in screening them for it.
The statistical test of interest then, is whether $\,p\,$ is moderate or whether it is extremely close to either 0 or 1.

In order to test this hypothesis we have decided to sequentially test patients until we have observed at least $\,c\,$ of both the wildtype (normal) and abnormal genotypes.
A small number of observations necessary to obtain at least $\,c\,$ of both genotypes is statistical evidence that $\,p\,$ is not far from 1/2.
Similarly, a large number of samples needed to observe at least $\,c\,$ of both genotypes is statistical evidence that the Bernoulli parameter $\,p\,$ is extreme.

Let $\,Y\,$ denote the `excess' number of trials needed beyond the minimum of $\,2c$.
The probability mass function of the maximum negative binomial distribution is
\begin{equation}                                                  
    \Pr[\,Y=y\,] = { {2c+y-1} \choose {c-1} } 
                (p^y + q^y) (pq)^c                            \label{MxNB}
\end{equation}
for $\,y=0,1,\ldots\,$ and $\,q=1-p$.

The maximum negative binomial distribution is so-named because it represents the larger of two negative binomial distributions: the number of failures before the $\,c-$th success is observed and the number of successes until the $\,c-$th failure is observed.
This distribution is also a mixture of two negative binomial distributions~(\ref{NB}) that are left-truncated at $\,Y=c\,$.

An intuitive description of the terms in~(\ref{MxNB}) are as follows.
There are $\,c\,$ successes and $\,c\,$ failures that occur with probability $\,(pq)^c\,$.
All of the $\,y\,$ extra trials beyond $\,2c\,$ must all be either successes or failures, hence the $\,p^y+q^y\,$ term.
Finally, the last $\,2c+y\,$ Bernoulli trial must be the one that completes the experiment ending with either the $\,c-$th success or the $\,c-$th failure.

In Zhang, {\em et al.} (2000) we describe properties of the distribution~(\ref{MxNB}).
The maximum negative hypergeometric distribution given at~(\ref{bino}) below and developed in the following sections is the finite sample analogue to the maximum negative binomial distribution~(\ref{MxNB}).

The parameters $\,p\,$ and $\,q=1-p\,$ are not identifiable in~(\ref{MxNB}).
Specifically, the same distribution in~(\ref{MxNB}) results when $\,p\,$ and $\,q\,$ are interchanged.
Similarly, it is impossible to distinguish between inference on $\,p\,$ and on $\,1-p\,$ without additional information.
In words, we can't tell if we are estimating $\,p\,$ or $\,q\,$ unless we also know how many successes and failures were observed at the point at which we obtained at least $\,c\,$ of each.
A similar identifiability problem is presented for the maximum negative hypergeometric distribution described in Section~4.

\begin{table}[t!]
\begin{center}
\caption{Comparison of related sampling distributions}
\medskip
\begin{tabular}{lll}
Sampling scheme \hspace*{.25in} & Infinite population & 
  Finite population \\[-1.5ex]
                & or with replacement \hspace*{.45in} & 
  without replacement \\ \hline
Predetermined   & Binomial distribution & Hypergeometric distribution
                                           \\[-1.5ex]
number of items                            \\[.5ex]
Until $\,c\,$   & Negative binomial    & Negative hypergeometric 
                                            \\[-1.5ex]
successes       & distribution (\ref{NB}) & distribution (\ref{NH})
                                            \\[.5ex]
Until either    & Riff shuffle or & Minimum negative   
                                           \\[-1.5ex]
$c\,$ successes or & Minimum negative  & hypergeometric distribution
                      (\ref{MinNH})        \\[-1.5ex]
$c\,$ failures  &  binomial distribution (\ref{riff}) &   
                                           \\[.5ex]
Until           &  Maximum negative & Maximum negative 
                                           \\[-1.5ex]
$\,c\,$ successes and  & binomial distribution (\ref{MxNB}) & 
    hypergeometric distribution (\ref{MxNH}), (\ref{bino})  
                                           \\[-1.5ex]
$\,c\,$ failures & &                       \\[.5ex] \hline
\end{tabular}
\end{center}
\end{table}


The {\em minimum negative binomial\/} or {\em riff-shuffle distribution\/} is the distribution of the smallest number of Bernoulli trials needed in order to observe either $\,c\,$ successes or $\,c\,$ failures.
Clearly, at least $\,c\,$ and fewer than $\,2c\,$ Bernoulli trials are necessary.
The random variable $\,Y+c\,$ counts the total number of trials needed until either $\,c\,$ successes or $\,c\,$ failures are observed for 
$\,Y=0,1,\ldots,c-1$.
The experiment ends with sample numbered $\,Y+c\,$ from the Bernoulli population.

The mass function of the minimum negative binomial distribution is
\begin{equation}                                                 
   \Pr[\,Y=y\,] = {{c+y-1} \choose {c-1}}
        \left( p^c q^y + p^y q^c \right)                     \label{riff}
\end{equation}
for $\,y=0,1,\ldots, c-1$.

The naming of~(\ref{riff}) as the minimum negative binomial refers to the smaller of two dependent negative binomial distributions: the number of failures before the $\,c-$th success, and the number of successes before the $\,c-$th failure.
In words, distribution~(\ref{riff}) says that there will be either $\,c\,$ Bernoulli successes and $\,y\,$ failures or else $\,c\,$ failures and $\,y\,$ Bernoulli successes.
This distribution is introduced by Uppuluri and Blot~(1970) and described in Johnson, Kotz, and Kemp (1992, pp 234--5).
Lingappaiah~(1987) discusses parameter estimation for distribution~(\ref{riff}).

The three discrete distributions described up to this point are based on sampling from an infinitely large Bernoulli $\,(p)\,$ parent population.
Each of these distributions also has a finite sample analogy.
These will be described next.

The {\em negative hypergeometric distribution\/} (Johnson, Kotz, and Kemp, 1992, pp 239--42) is the distribution of the number of unsuccessful draws from an urn with two different colored balls until a specified number of successful draws have been obtained.
If $\,m\,$ out of $\,N\,$ balls are of the `successful' type then the number of unsuccessful draws $\,Y\,$ observed before $\,c\,$ of the successful types are obtained is
\begin{equation}                                                  
  \Pr[\,Y=y\,] = { {c+y-1} \choose {c-1} } 
                 { {N-c-y} \choose {m-c} }
           \left/ { {N} \choose {m} } \right.                   \label{NH}
\end{equation}
with parameters satisfying $\,1 \leq c \leq m <N\,$ and range $\,y=0,1,\ldots, N-m$.
The expected value of $\,Y\,$ in~(\ref{NH}) is $\,mc/(N-m-1)\,$.

The negative hypergeometric distribution~(\ref{NH}) is the finite sample analogy to the negative binomial distribution~(\ref{NB}).
Unlike the negative binomial distribution, the negative hypergeometric distribution has a finite range.
The maximum negative hypergeometric distribution described in the following sections is the larger of two, dependent negative hypergeometric distributions.

The {\em minimum negative hypergeometric distribution\/} describes the smallest number of urn draws needed in order to observe either $\,c\,$ successes or $\,c\,$ failures.
This distribution is the finite sample analogy to the riff-shuffle distribution~(\ref{riff}).
The probability mass function of the minimum negative hypergeometric distribution is
\begin{equation}                                               
 \Pr[\,Y=y\,] = {{c+y-1}\choose{c-1}}\!\!
     \left\{\!
        {{m}\choose{c}} \! {{N-m}\choose{y}} +
        {{m}\choose{y}} \! {{N-m}\choose{c}}\!\right\}       \label{MinNH}
     \!\!\left/\!
     \left\{\!
        {{c+y}\choose{c}}\! {{N}\choose{c+y}}
     \!\right\} 
     \right.
\end{equation}
for $\,y=0,1,\ldots, c-1$.

In the example of the charge to Noah, we have $\,c=7\,$ male/female pairs of animals captured from a finite population of $\,m\,$ males and $\,N-m\,$ females.

In Section~2 we give the probability mass function of the maximum negative hypergeometric distribution.
Section~3 details some approximations to this distribution.
In Section~4 we discuss estimation of the parameter that describes the contents of the urn.


\section                  {The distribution}

An urn contains $\,N\,$ balls: $\,m\,$ of one color; and the 
remaining $\,N-m\,$ of another color.
We continue sampling from the urn without replacement until we have observed $\,c\,$ balls of both colors, for integer parameter $\,c\geq 1$.
Sampling with replacement is the same as sampling from the maximum negative binomial distribution~(\ref{MxNB}) with parameter $\,p=m/N\,$

Let $\,Y\,$ denote the random variable counting the number
of extra draws needed beyond the minimum $\,2c\,$.
That is, on draw numbered $\,Y+2c\,$ we will have first observed at least $\,c\,$ of both colors.
All of the $\,Y\,$ extra draws from the urn must be of the same color so there will be $\,c\,$ of one color and $\,Y+c\,$ of the other color at the end of the experiment.
We will describe the distribution and properties of this random variable.

For $\,k=1,2,\ldots\,$ define the {\em factorial polynomial}
$$
      z^{(k)} = z(z-1)\cdots(z-k+1)\; .
$$
We also define $\,z^{(0)}=1$.

The {\em maximum negative hypergeometric distribution \/} probability mass function can be written as
\begin{equation}                                                 
  \Pr[\,Y=y\,] = \left.{{2c+y-1}\choose{c-1}}
        \left\{m^{(c+y)}(N-m)^{(c)} 
            + m^{(c)}(N-m)^{(c+y)}\right\}                    \label{MxNH}
        \right/ N^{(2c+y)}
\end{equation}
defined for the range of $\,Y$:
$$
   0  \leq  y  \leq  \max \{ m-c,\, N-m-c\} \; .
$$

The integer valued parameters $\,(N,m,c)\,$ are constrained to
$$
     1 \leq c \leq m < N {\rm \ \ and \ \ } c \leq N-m \; .
$$

Similarly,
\begin{equation}                                                 
 \Pr[\,Y=y\,] =  \{ c/(2c+y)\} \left\{
           {{m}\choose{c+y}}{{N-m}\choose{c}}
         + {{m}\choose{c}}{{N-m}\choose{c+y}}\right\}
       \left/ {{N}\choose{2c+y}}\right.                     \label{bino}
\end{equation}
expresses the maximum negative hypergeometric distribution~(\ref{MxNH}) in terms of binomial co\-efficients.

The same distribution in~(\ref{MxNH}) and~(\ref{bino}) result when the parameter $\,m\,$ is interchanged with $\,N-m$.
This remark illustrates the identifiability problem with the parameters in the maximum negative hypergeometric distribution.
A similar identifiability problem occurs in the maximum negative binomial distribution given at~(\ref{MxNB}).
We will describe the estimation of the $\,m\,$ parameter in Section~4.

Special cases of this distribution are as follows.
For general parameter values,
$$
  \Pr[\,Y=0\,] = \left.{{N-2c}\choose{m-c}} {{2c}\choose{c}}
                 \right/ {{N}\choose{m}} \; .
$$

If $\,c=m=N/2\,$ then the maximum negative hypergeometric distribution is degenerate and all of its probability is a point mass at $\,Y=0$.
In words, if $\,c=m=N-m\,$ then there can be only one possible outcome.
In this case, all of the balls in the urn must be drawn before we can observe $\,c\,$ balls of both colors.

The special case of $\,c=m=1\,$ with $\,N>2\,$ has the form
$$
 \Pr[\,Y=y\mid m=c=1\,] = \left\{
     \begin{tabular}{ll}
          $2/N$ & for $\,y=0$ \\
          $1/N$ & for $\,y=1,\ldots ,N-2$
     \end{tabular}    \right.
$$
and zero otherwise.
This is also the form of the distribution for $\,c=1\,$ and $\,m=N-1$.

The special case for $\,c=m\,$ and $\,N=2m+1\,$ has mass function
$$
 \Pr[\,Y=y\mid m=c,\, N=2c+1\,] = \left\{
     \begin{tabular}{rl}
          $(m+1)/(2m+1)$ & for $\,y=0$ \\
          $    m/(2m+1)$ & for $\,y=1$
     \end{tabular}    \right.
$$
and zero otherwise.
This is also the distribution of $\,Y\,$ for $\,m=c+1\,$ and $\,N=2m+1$.
In words, this represents the distribution of the color of the last ball remaining after all but one have been drawn from the urn.


\section             {Properties and Approximations}


\begin{table}[t!]
\caption{Ranges of $\,m\,$ parameter that result in unimodal maximum negative hypergeometric distributions for specified values of $\,c\,$ and $\,N\,$.
Omitted distributions are either degenerate or the parameter values are invalid }
                                                       \label{modetable}
\begin{center}
\begin{tabular}{cccc}
\makebox[.75in]{$c$} & \makebox[1.25in]{$N=10$} & \makebox[1.25in]{$N=50$}  & \makebox[1.25in]{$N=250$} \\ \hline
\p 1 & \multicolumn{3}{c}{Unimodal for all $\,m=1,\ldots ,N-1$}\\
\p 2 & $3\leq m\leq 7$ & $ 9\leq m\leq 41$ & $ 38\leq m\leq 212$ \\
\p 3 & $4\leq m\leq 6$ & $13\leq m\leq 37$ & $ 55\leq m\leq 195$ \\
\p 4 & $m=5$           & $15\leq m\leq 35$ & $ 65\leq m\leq 185$ \\
\p 5 & $---$           & $16\leq m\leq 34$ & $ 73\leq m\leq 177$ \\
  10 & $---$           & $20\leq m\leq 30$ & $ 90\leq m\leq 160$ \\
  15 & $---$           & $22\leq m\leq 28$ & $ 98\leq m\leq 152$ \\
  20 & $---$           & $24\leq m\leq 26$ & $103\leq m\leq 147$ \\
  25 & $---$           & $---$             & $106\leq m\leq 144$ \\ \hline
\\
\end{tabular}
\end{center}
\end{table}


There are five basic shapes that the maximum negative hypergeometric distribution will assume.
These are illustrated in Figs.~1 through~5.
In each figure, the limiting maximum negative binomial distribution~(\ref{MxNB}) is also presented.
This limit can be expressed, more formally,  as follows.

{\bf Lemma 1.}
{\em For fixed values of $\,c\geq 1$, let $\,N\,$ and $\,m\,$ both grow large such that $\,m/N=p\,$ for $\,p\,$ bounded between zero and one.
Then the behavior  of the maximum negative hypergeometric random variable~{\rm(\ref{MxNH})} approaches the maximum negative binomial distribution~{\rm(\ref{MxNB})} with parameters c and p.}

{\bf Proof.}
Values of $\,Y\,$ remain bounded with high probability under these conditions.
In~(\ref{MxNH}) we write
\begin{eqnarray*}
 m^{(c+y)} \, (N-m)^{(c)} \left/ N^{(2c+y)}\right. & \geq &
   (m-c-y)^{c+y} \, (N-m-c)^c \left/ N^{2c+y} \right.  \\
    &\hspace*{-1.5in}=& \hspace*{-.75in}(m/N)^{c+y}\, \{(N-m)/N\}^c 
         \{1-(c+y)/m\}^{c+y} \{1-c/(N-m)\}^c     \\
    &\hspace*{-1.5in}=& \hspace*{-.75in}p^{c+y}\, q^c \{1+O_p(N^{-1})\} 
\end{eqnarray*}
where $\,p=m/N\,$ and $\,q=1-p=(N-m)/N$.

We can also write
\begin{eqnarray*}
 m^{(c+y)}\, (N-m)^{(c)} \left/ N^{(2c+y)} \right. & \leq &
   m^{c+y} \, (N-m)^c \left/ (N-2c-y)^{2c+y}\right. \\
  &=& (m/N)^{c+y}\,\{(N-m)/N\}^c \left/ \{1-(2c+y)/N\}^{2c+y}\right. \\
 &=& p^{c+y}\, q^c \{1+O_p(N^{-1})\} \; .
\end{eqnarray*}

A similar argument shows 
$$
 m^{(c)}\,(N-m)^{(c+y)}\left/ N^{(2c+y)} \right. =
    p^c\,q^{c+y} \{1+O_p(N^{-1})\} 
$$
completing the proof.~\rule{.12in}{.12in}

In words, if $\,m\,$ and $\,N\,$ are both large then sampling from the urn without replacement is almost the same as sampling with replacement.
Sampling with replacement is the same as sampling from a Bernoulli parent population yielding the maximum negative binomial distribution~(\ref{MxNB}).

We next describe the modes for this distribution.
The maximum negative hypergeometric distribution can have either one or two modes.
Write
$$
   \Pr[\,Y=0\,] / \Pr[\,Y=1\,] = (c+1)/c
$$
to show that this distribution always has at least one local mode at $\,Y=0$.

The maximum negative binomial distribution~(\ref{MxNB}) also has at least one local mode at $\,Y=0\,$ for all values of the parameter $\,p$.
The local mode of the maximum negative hypergeometric distribution at $\,Y=0\,$ is clearly visible in Figs. 1, 2, 4, and 5.
The local mode at $\,Y=0\,$ in Fig.~3 is also present but it is very small.

Table~\ref{modetable} presents examples of parameter values corresponding to unimodal distributions in~(\ref{bino}).
In general, there will be only one mode at $\,Y=0\,$ when $\,m/N\,$ is not too far from 1/2.
The range of $\,m\,$ with unimodal distributions becomes narrower as $\,c\,$ becomes larger when $\,N\,$ is fixed.
If $\,m=N/2\,$ then the distribution is always unimodal.


\subsection            {A gamma approximation}

An approximate gamma distribution is illustrated in Fig.~4.
Under the conditions of the following lemma, the local mode at $\,Y=0\,$ becomes negligible.

{\bf Lemma 2.}
{\em For fixed $\,c\geq 1,$ if $\,m\,$ grows as $\,\theta N^{1/2}\,$ for large $\,N\,$ and some $\,\theta>0\,$ then $\,\theta Y/N^{1/2}\,$ behaves approximately as the sum of $\,c\,$ independent standard exponential random variables.}

{\bf Proof.}
Begin at~(\ref{MxNH}) and write
\begin{eqnarray*}
   {{2c+y-1} \choose {c-1}} &=& \prod_{i=1}^{c-1} (y+2c-i)/i \\
    &=& y^{c-1}(1+\BO_p(N^{-1/2}))/\Gamma(c)
\end{eqnarray*}

Define $\,\Delta\,$ as
$$
   \Delta = m^{(c)} (N-m)^{(y+c)}/N^{(2c+y)} \; .
$$

Under the conditions of his lemma, the term
$$
         m^{(c+y)} (N-m)^{(c)}/N^{(2c+y)}
$$
will be much smaller than $\,\Delta\,$ and can be ignored.

We have
\begin{eqnarray*}
  \log \Delta &=& \sum_{i=0}^{c-1}  \,\log\{(m-i)/(N-i)\}
        + \sum_{j=0}^{y+c-1} \,\log\{(N-m-j)/(N-c-j)\} \\
  &=& c \log\{\theta N^{-1/2} + \BO(N^{-1})\} 
         + \sum_{j=0}^{y+c-1} \,\log\{ 1-(m-c)/(N-c-j)\}
\end{eqnarray*}

For $\,\epsilon\,$ near zero, write
$$
  \epsilon - \epsilon^2/2 \,\leq\, \log(1+\epsilon) 
         \,\leq\,\epsilon
$$
so that
$$
  \log\Delta = c\log\{\theta N^{-1/2}+ \BO(N^{-1})\} 
               -\theta y /N^{1/2} + \BO_p(y/N) \;
$$

The transformation $\,X=\theta Y/N^{1/2}\,$ has Jacobian $\,N^{1/2}/ \theta\,$ so
$$
 (N^{1/2}/\theta) \Pr[\,Y\,] \; = \;
      x^{c-1} {\rm e}^{-x}/\Gamma(c)
$$
ignoring terms that tend to zero for large vales of $\,N.$
This is the density function of the sum of $\,c\,$ independent, standard exponential random variables.~\rule{.12in}{.12in}


\subsection             {A half-normal approximation}

If $\,Z\,$ has a standard normal distribution then the distribution of $\,\mid Z \mid\,$ is said to be {\em standard half-normal} or {\em folded normal}.
The density function of the random variable $\,X=\mid Z \mid\,$ is
$$
  (2/\pi)^{1/2}\exp(-x^2/2)
$$
for $\,x\geq 0\,$ (Stuart and Ord, 1987, p 117).
The approximate half-normal behavior of the maximum negative hypergeometric distribution is illustrated in Fig.~5.

{\bf Lemma 3. } {\em When $\,N\,$ becomes large, if $\,m=N/2\,$ and $\,c\,$ grows as $\,N^{1/2}\,$ then $\,Y/(2c)^{1/2}\,$ behaves approximately as a standard half-normal random variable.}

The proof involves expanding all factorials in~(\ref{bino}) using Stirling's approximation.
The details are provided in Appendix A.


\subsection             {A normal approximation}

The normal approximation to the maximum negative hypergeometric distribution can be seen in Fig.~3.
This is proved more formally in Lemma~4, below.
No generality is lost by requiring $\,m>N/2\,$ because $\,m\,$ and 
$\,N-m\,$ can be interchanged to yield the same distribution.

{\bf Lemma 4. } {\em For large values of $\,N$, suppose $\,c\,$ grows as $\,N^{1/2}\,$ and $\,m=Np\,$ for $\,1/2 < p < 1.$
Then $\,(Y-\mu)/\sigma\,$ behaves approximately as standard normal where
$$
  \mu = c(p-q)/q
$$
for $\,q=1-p\,$ and
$$  
   \sigma = (cp)^{1/2}/q \; .
$$
}

The proof of thsi lemma is given in Appendix B.
The details involve using Stirling's approximation to all of the factorials in~(\ref{bino}) and expanding these in a two-term Taylor series.


\section                    {Estimation}

The most practical situation concerning parameter estimation involves estimating the $\,m\,$ parameter when $\,c\,$ and $\,N\,$ are both known.
In terms of the original, motivating example drawing inference on the genetic markers in cancer patients, the finite population size $\,N\,$ will be known, and the parameter $\,c\,$ is chosen by the investigators in order to achieve specified power and significance levels.
The $\,m\,$ parameter describes the composition of the $\,N\,$ individuals in the finite-sized population.
The value of $\,m\,$ is known without error if all $\,N\,$ subjects are observed.

The estimation of $\,m\,$ in this section is made on the basis of a single observation of the random variable $\,Y$.
We will treat the unknown $\,m\,$ parameter as continuous valued rather than as a discrete integer as it has been used in previous sections.

The log-likelihood kernel function of $\,m\,$ in~(\ref{MxNH}) is
$$
  \Lambda(m) = \log\left\{\left(m^{(c)}(N-m)^{(c+y)} +
             m^{(c+y)}(N-m)^{(c)} \right)
           \left/ N^{(2c+y)} \right. \right\}  \; .
$$

As a numerical illustration, the function $\,\Lambda(m)\,$ is plotted in Fig.~6 for $\,N=20\,$ and $\,c=3$.
Observed values of $\,y\,$ are given as $\,0,1,\ldots, 7\,$ in this figure.
The range of valid values of the $\,m\,$ parameter are $\,3\leq m\leq 17\,$ for the values of $\,c\,$ and $\,N\,$ in this example.
Smaller observed values of $\,y=0,1,2\,$ in this example exhibit log-likelihood functions with a single mode corresponding to maximum likelihood estimates of $\,\widehat{m}=N/2\,$
For values of $\,y\geq 3\,$ the likelihood $\,\Lambda\,$ has two modes, symmetric about $\,N/2\,$.

Intuitively, if the observed value $\,Y\,$ is small then we are inclined to believe that the urn is composed of an equal number of balls of both colors.
That is, if we quickly observe $\,c\,$ of both colored balls then this is good statistical evidence of an even balance of the two colors in the urn.
Conversely, if the observed $\,Y\,$ is relatively large then we will estimate an imbalance in the composition of the urn.
Without the additional knowledge of the number of successes and failures observed then we are unable to tell if we are estimating $\,m\,$ or $\,N-m$.

More generally, there will be either one mode of $\,\Lambda(m)\,$ at $\,\widehat{m}=N/2\,$ or else two modes, symmetric about $\,N/2\,$ depending on the sign of
$$
   \Lambda''(m) = (\partial/\partial m)^2 \, \Lambda(m)
$$
evaluated at $\,m=N/2$.
If $\,\Lambda''(N/2)\,$ is negative then there will be one mode of $\,\Lambda\,$ at $\,N/2$.

Useful rules for differentiating factorial polynomials are as follows.
For $\,c\geq 1$,
$$
   (\partial/\partial m) \; m^{(c)} =
       m^{(c)} \,\sum_{i=0}^{c-1}\, (m-i)^{-1}
$$
and for $\,c\geq 2$,
$$
    (\partial/\partial m)^2 \; m^{(c)} =   2m^{(c)}\,
       \sum_{\stackrel{\scriptstyle{i=0}}{i<i'}}^{c-1}\, 
         (m-i)^{-1} (m-i')^{-1} \; .
$$

Use these rules to write
$$
  (\partial /\partial m) \; m^{(c)} (N-m)^{(y+c)} =
       m^{(c)} (N-m)^{(y+c)} \left\{ \,\sum_{i=0}^{c-1}\, (m-i)^{-1} 
      \; - \;\sum_{j=0}^{y+c-1}\, (N-m-j)^{-1}
     \right\}
$$
and show that the likelihood $\,\Lambda(m)\,$ always has a critical value at $\,m=N/2$:
$$
  (\partial /\partial m)\Lambda(m)\Big|_{m=N/2} = 0 \; .
$$

The critical point of $\,\Lambda\,$ at $\,m=N/2\,$ may either be a global maximum or else a local minimum as seen in the example of Fig.~6.
This distinction depends on the sign of the second derivative $\,\Lambda''\,$ of $\,\Lambda\,$.

The second derivative of $\,\Lambda\,$ can be found from
\begin{eqnarray*}
   (\partial /\partial m)^2 \; m^{(c)} (N-m)^{(y+c)} & = &
        m^{(c)} (N-m)^{(y+c)} \\
 && \hspace*{-2in}\times\left[\left\{\,
   \sum_{i=0}^{c-1}\, (m-i)^{-1} 
    - \,\sum_{j=0}^{y+c-1}\, (N-m-j)^{-1} \right\}^2
    - \sum_{i=0}^{c-1}\, (m-i)^{-2} 
    - \,\sum_{j=0}^{y+c-1}\, (N-m-j)^{-2}   \right] \; .
\end{eqnarray*}

The sign of $\,\Lambda''(N/2)\,$ is the same as that of
$$
  \phi(N,c,y) =
       \sum_{\stackrel{\scriptstyle{k=0}}{k<k'}}^{y-1}\, 
         (N/2-c-k)^{-1}   (N/2-c-k')^{-1}
     \; - \;\sum_{i=0}^{c-1}\, (N/2-i)^{-2} \; .
$$

The first summation in $\,\phi\,$ is zero when $\,y=0\,$ or 1.
The function $\,\phi\,$ is then negative for small values of $\,y\,$ demonstrating that the maximum likelihood estimate of $\,m\,$ is $\,\widehat{m}=N/2\,$ in these cases.
Similarly, $\,\phi\,$ is an increasing function of $\,y\,$ and may eventually become positive for larger values of $\,y,$ so that $\,\Lambda\,$ will have two modes.
These modes are symmetric about $\,N/2\,$ because 
$$
    \Lambda(m) = \Lambda(N-m)
$$
for all $\,0<m<N$.

In other words, a small observed value of $\,y\,$ leads us to believe that there are an equal number of balls of both colors in the urn and estimate $\,m\,$ by $\,N/2$.
Similarly, a large observed value of $\,y\,$ relative to $\,c\,$ leads us to estimate an imbalance in the composition of the urn.


\section*                             {\bf References}

\begin{description}

\item Johnson, N.L., S. Kotz, and A.W. Kemp (1992). 
{\it Univariate Discrete Distributions\/}. New York: 
John Wiley \& Sons.

\item Lingappiah, G.S. (1987). Some variants of the binomial
distribution. {\it Bulletin of the Malaysian Mathematical 
Society\/} {\bf 10}: 82--94.

\item Stuart, A. and J.K. Ord (1987). {\em Kendall's Advanced Theory of Statistics\/} Vol 1, 5th Edition: {\em Distribution Theory,\/} New York: Oxford University Press.

\item Uppuluri, V.R.R., and W.J. Blot (1970). ``A probability distribution arising in a riff-shuffle.'' {\it Random Counts  in Scientific Work, 1: Random Counts in Models and Structures},  G.P. Patil (editor), University Park: Pennsylvania State University Press, pp  23--46.

\item Zhang, Z., B.A.~Burtness, and D.~Zelterman (2000). The maximum negative binomial distribution. {\em Journal of Statistical Planning and Inference\/} {\bf 87:} 1--19.

\end{description}


\section*{Appendix A: Proof of the half-normal approximation}

The proof of Lemma~3 is provided here.
We assume that $\,m=N/2\,$ and $\,c=N^{1/2}.$
The random variable $\,Y\,$ is $\,\BO_p(N^{1/4}).$

Expand all of the factorials in~(\ref{bino}) using Stirling's approximation giving
$$
 \log \Pr[\,Y=y\,] = -1/2\log(2\pi) + T_1(c) + T_2(N) + \BO_p(N^{-1/2})
$$
where
\begin{eqnarray*}
  T_1 & = & \log\{2c/(2c+y)\} - (c+y+1/2)\log(c+y) \\
      &   & -\, (c+1/2)\log(c) + (2c+y+1/2)\log(2c+y) 
\end{eqnarray*}
contains terms in $\,\BO_p(c\log c)\,$ and
\begin{eqnarray*}
  T_2 &=&  (N+1)\log(N/2) - (N+1/2)\log (N)
                 - (N/2 -c+1/2)\log(N/2-c) \\
      & &  -\, (N/2-c-y+1/2)\log(N/2-c-y) + (N-2c-y+1/2)\log(N-2c-y)
\end{eqnarray*}
contains terms that are $\,\BO_p(N\log N)\,$.

In all of the following expansions it is useful to keep in mind that $\,c=c_N\,$ is approximately equal to $\,N^{1/2}\,$ and $\,Y=\BO_p(N^{1/4})\,$.
Write $\,T_1\,$ as
\begin{eqnarray*}
 T_1 &=& -\log\{(2c+y)/2c\} + (2c+y+1/2)\log 2 -1/2\log(c) \\
     & &  \quad\quad +\, c\log\{(c+y/2)^2/c(c+y)\} + 
             (y+1/2)\log\{(c+y/2)/(c+y)\} \\
     &=& -\log(1+y/2c) + (2c+y+1/2)\log2 - 1/2\log(c) \\
     & &  \quad\quad +\, c\log\{1+y^2/4c(c+y)\} + 
             (y+1/2)\log\{1-y/2(c+y)\} \; .
\end{eqnarray*}

Expand every appearance of $\,\log(1+\epsilon)=\epsilon+\BO(\epsilon^2)\,$ for $\,\epsilon\,$ near zero to show
$$
   T_1 = (2c+y+1/2)\log 2 - 1/2\log(c) -y^2/4c + \BO_p(N^{-1/4})\; .
$$

Similarly, we can write
\begin{eqnarray*}
 T_2 &=&   (N/2-c)\log\{(N-2c-y)/(N-2c)\} \\
     & &  \quad\quad +\, y\log\{(N-2c-2y)/(N-2c-y)\} -(2c+y-1)\log2 \\
     & &  \quad\quad +\, 1/2\,\log\{N/(N-2c)\} + 
             1/2\,\log\{(N-2c-y)/(N-2c-2y)\} \\
     &=& (N/2-c)\log\{1+y^2/\BO_p(N^2)\} + y\log\{1-y/\BO_p(N)\} \\
     & & \quad\quad -\, (2c+y-1)\log2 + 1/2\,\log\{1+\BO(c/N)\} 
            + 1/2\,\log\{1+y/\BO_p(N)\} \; .
\end{eqnarray*}

Then write $\,\log(1+\epsilon)=\BO(\epsilon)\,$ for $\,\epsilon\,$ near zero, giving
$$
   T_2 = -(2c+y-1)\log2 + \BO_p(N^{-1/2}) \; .
$$

These expressions for $\,T_1\,$ and $\,T_2\,$ in $\,\Pr[Y]\,$ give
$$
   \log \Pr[\,Y=y\,] = 
         \log 2 -1/2\log (\pi c)- y^2/4c + \BO_p(N^{-1/4})\; . 
$$

Finally, we note that $\,(2c)^{1/2}\,$ is the Jacobian of the transformation $\,X=Y/(2c)^{1/2}\,$.
Then
$$
   \log \{ (2c)^{1/2}\Pr[\,Y=y\,]\} = 
      1/2\log (2/\pi) - x^2/2 + \BO_p(N^{-1/4})
$$
is the density of the folded normal distribution, except for terms that tend to zero with high probability.~\rule{.12in}{.12in}


\section*{Appendix B: Standard normal approximate distribution}

The details of the proof of Lemma 4 are given here.
Define $\,\Omega\,$ as
$$
   \Omega = \{c/(2c+y)\} {{m}\choose{c+y}}
{{N-m}\choose{c}} \left/ {{N}\choose{2c+y}}\right. \; .
$$

The term
$$
  \{c/(2c+y)\} {{m}\choose{c}}
{{N-m}\choose{c+y}} \left/ {{N}\choose{2c+y}}\right. \; .
$$
in~(\ref{bino}) is much smaller than $\,\Omega\,$ and can be ignored under the conditions of this lemma.

Expand all of the factorials in $\,\Omega\,$ using Stirling's formula giving
$$
    \log \Omega = - 1/2 \log(2\pi) + S_1(N) + S_2(N-c) + S_3(y) 
        + \BO_p(N^{-1/2})
$$
where
$$
   S_1 = (Np+1/2)\,\log(Np) + (Nq +1/2)\,\log(Nq)
           - (N+1/2)\,\log N
$$
corresponds to $\,m!\, , (N-m)!\, , $ and $\,N!\,$;
\begin{eqnarray*}
   S_2 & = & (N-2c-y+1/2)\,\log(N-2c-y)
             - (Np-c-y+1/2)\,\log(Np-c-y) \\
       && \quad  - (Nq -c+1/2)\,\log(Nq-c)
\end{eqnarray*}
corresponds to $\,(N-2c-y)!\, , (m-c-y)!\, , $ and $\,(N-m-c)!\,$; and
\begin{eqnarray*}
   S_3 & = & \log\{c/(2c+y)\} + (2c+y+1/2)\,\log(2c+y) \\
   &&\quad -(c+y+1/2)\,\log(c+y) - (c+1/2)\,\log(c)
\end{eqnarray*}
corresponds to $\,c/(2c+y),\, (2c+y)!,\, (c+y)!\, $ and $\,c!\,$ .

Write out all of the terms in $\,S_1\,$ to show
$$
   S_1 = Np\,\log p + Nq\,\log q + 1/2 \log(Npq) \; .
$$

We can write $\,S_2\,$ as
\begin{eqnarray*}
  S_2 &=& (N-2c-y+1/2)\,\log N  + (N-2c-y+1/2)\,\log\{1-(2c+y)/N)\} \\
&&\quad -(Np-c-y+1/2)\,\log (Np) -(N-c-y+1/2)\,\log\{1-(c+y)/Np)\} \\
&&\quad -(Nq-c+1/2)\,\log (Nq) -(Nq-c+1/2)\,\log(1-c/N) \; .
\end{eqnarray*}

Then
\begin{eqnarray*}
 S_1 + S_2 &=& c\log(pq) + y\log p 
       + (N-2c-y+1/2)\,\log\{1-(2c+y)/N\} \\
  &&\quad -(Np-c-y+1/2)\,\log\{1-(c+y)/Np\}
        -(Nq-c+1/2)\,\log(1-c/Nq) \; .
\end{eqnarray*}

Since
$$
   (c+y)/N = \BO_p(N^{-1/2})
$$
we can expand
$$
    \log(1+\epsilon) = \epsilon - \epsilon^2/2 + \BO(\epsilon^3)
$$
for $\,\epsilon\,$ near zero to show
\begin{eqnarray*}
  S_1 + S_2 &=& c\log (pq) + y \log p + (2c+y)^2/2N \\
    && \quad -(c+y)^2/2Np - c^2/2Nq + \BO_p(N^{-1/2}) \; . 
\end{eqnarray*}

Then write
$$
                y=\mu + Z \sigma
$$
where $\,Z\,$ is a $\,\BO_p(1)\,$ random variable, giving
\begin{eqnarray*}
  S_1 + S_2 &=& c\log(pq) + y\log p + 
            \{ (2c+\mu)^2 pq - (c+\mu)^2q - c^2p\}/2Npq \\
   && \quad + Z\sigma\{(2c+\mu)p - c - \mu\}/Np 
            + \BO_p(N^{-1/2})\; .
\end{eqnarray*}

Substitute $\,\mu=c(p-q)/q\,$ to show
\begin{equation}                                                  
  S_1 + S_2 = c\log(pq) + y \log p + \BO_p(N^{-1/2})  \; .    \label{1+2}
\end{equation}

Next write $\,S_3\,$ as
$$
  S_3 = (c+y)\log\{(2c+y)/(c+y)\} + (c-1/2)\log\{(2c+y)/c)\}
     -1/2 \log(c+y) \; .
$$

Expand the argument of the first logarithm in $\,S_3\,$ here in a two-term Taylor series showing
\begin{eqnarray*}
 (2c+y)/(c+y) &=& (2c + \mu + Z\sigma)/(c+\mu+Z\sigma)\\
 &\hspace*{-1in}=& \hspace*{-.5in} \{(2c+\mu)/(c+\mu)\} 
      \; [ \{1+Z\sigma/(2c+\mu)\} \, / \, 
                   \{1+Z\sigma/(c+\mu)\} ]\\
      & \hspace*{-1in}=& \hspace*{-.5in} \{(2c+\mu)/(c+\mu)\} \\
      &  \hspace*{-.95in} & \hspace*{-.5in} \times\;
          \left[1 \; - \;Z\sigma c/\{(c+\mu)(2c+\mu)\}\; + \;
         Z^2\sigma^2c/\{(c+\mu)^2(2c+\mu)\} + \BO(N^{-3/4}) \right] \; . 
\end{eqnarray*}

Then expand 
$$
    \log(1+\epsilon)= \epsilon -\epsilon^2/2 + \BO(\epsilon^3)
$$
to show
\begin{eqnarray*}
(c+y)\log\{(2c+y)/(c+y)\} &=& (c+\mu+Z\sigma)\log\{(2c+\mu)/(c+\mu)\} \\
   && -\; (c+\mu+Z\sigma)\,Z\sigma c/\{(c+\mu)(2c+\mu)\} \\
   && +\; Z^2\sigma^2c/\{(c+\mu)(2c+\mu)\}  \\
   && -\;1/2 \, (Z\sigma c)^2/\{2(c+\mu)(2c+\mu)^2\} + \BO_p(N^{-1/4}) \; .
\end{eqnarray*}

Similarly,
\begin{eqnarray*}
 (c-1/2)\log\{(2c+y)/c)\} &=& 
      (c-1/2)\log (2+\mu/c) + (c-1/2)\log\{1+ Z\sigma/(2c+\mu)\} \\
  &=&(c-1/2)\log(2+\mu/c) +  Z\sigma c/(2c+\mu) \\
  && \quad\quad - \; Z^2\sigma^2 c /\{2(2c+\mu)^2\} + \BO_p(N^{-1/4})
\end{eqnarray*}
and
\begin{eqnarray*}
 1/2\, \log (c+y) &=& 1/2 \log\{(c+\mu)(1+Z\sigma/(c+\mu))\} \\
    &=& 1/2 \log(c+\mu) + \BO_p(N^{-1/4})
\end{eqnarray*}
so that
\begin{eqnarray*}
 S_3 &=& (c+y)\,\log\{(2c+\mu)/(c+\mu)\} - 1/2\log(c+\mu) \\
   & & +\;(c-1/2)\log\{(2c+\mu)/c\} -Z^2\sigma^2c/\{2(c+\mu)(2c+\mu)\}
        + \BO_p(N^{-1/4})\; .
\end{eqnarray*}

Substitute values of $\,\mu=c(p-q)/q\,$ and $\,\sigma^2=cp/q^2\,$ giving
$$
   S_3 = (c+y)\log p - \log\sigma + (c-1/2)\log q 
          - Z^2/2  + \BO_p(N^{-1/4}) \; .
$$

This expression, together with the form of $\,S_1+S_2\,$ given at~(\ref{1+2}) shows that
$$
\log \Omega = -1/2\log(2\pi) -\log\sigma -Z^2/2 + \BO_p(N^{-1/4})
$$
demonstrating the approximately standard normal behavior of the random variable $\,Z=(Y-\mu)/\sigma\,$.~\rule{.12in}{.12in}




\begin{figure}[ht]
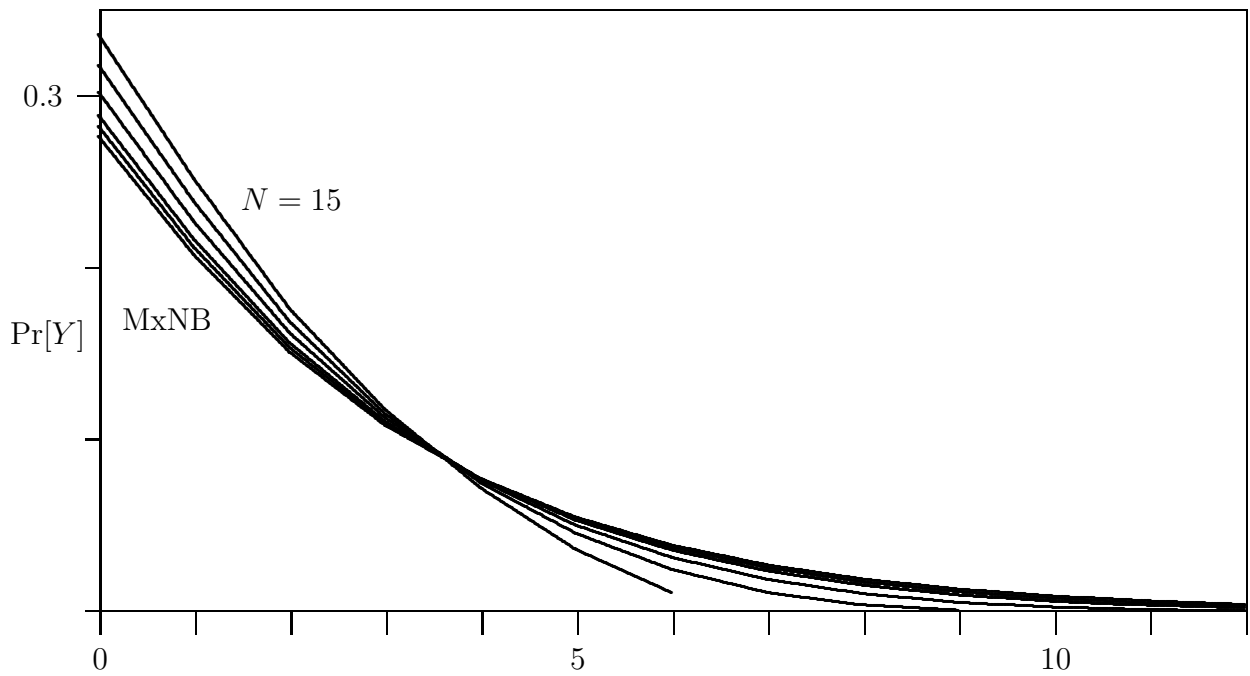


\caption{The maximum negative hypergeometric distribution 
for $\,c=3\,$ and $\,m=0.4N\,$ with values of $\,N=15,\; 20,\; 
30,\; 60,\,$ and 120. The distribution corresponding to the maximum negative binomial (MxNB) given at~(\ref{MxNB}) with parameters $\,c=3\,$ and $\,p=0.4\,$ is the limit when $\,N\,$ is large.  Mass points of these discrete distributions are joined with lines for clarity in this and other figures. }
\typeout{ Figure 1}
\bigskip

\beginpicture
\setcoordinatesystem units < .5in, 9in>
\setplotarea x from 0 to 12, y from 0  to  .35

\axis top   /
\axis right /

\axis bottom ticks out
    numbered from 0 to 10 by 5  
    unlabeled from 0 to 12 by 1  / 

\axis left ticks out
    numbered from  .3 to  .3 by  .3
    unlabeled short from 0 to .2 by .1 /

\put{$N=15$} at 2  .24
\put{MxNB}    at .7 .17
\put {$\Pr[Y]$} at -.55 .16

\setlinear

\plot          
0 .33566 1 .2517482 2 .17622 3 .11748 4 .07192 5 .03596 6 .01098 /

\plot    
0 .317853 1 .238390 2 .168706 3 .114916 4 .074874
5 .045367 6 .024450 7 .01100  8 .003715 9 .000722 /

\plot
 0 .3023367437 
 1 .2267525578  2 .1621773729  3 .1127190251  4 .0763609404
 5 .0500782576  6 .0314357570  7 .0186539003  8 .0103310751
 9 .0052668403 10 .0024303638 11 .0009919852 12 .0003459230 /

\plot     
 0 .288658525  1 .216493894  2 .156311315 3 .110561662 4 .077062683
 5 .052920524  6 .035709493  7 .023603558 8 .015242620 9 .009597695
10 .005883938 11 .003508118 12 .002032146 /

\plot    
 0 .282399486  1 .211799614  2 .153576917 3 .109474193 4 .077204994
 5 .053930967  6 .037279667  7 .025463284 8 .017164604 9 .011410087
10 .007476804 11 .004829270 12 .003074922 /

\plot 
 0 .276480000  1 .207360000  2 .150958080 3 .108380160 4 .077237453
 5 .054743040  6 .038587761  7 .027035763 8 .018817485 9 .013007353
10 .008928989 11 .006087946 12 .004124032 /

\endpicture
\end{figure}


\typeout{ Figure 2}
\begin{figure}[hb]
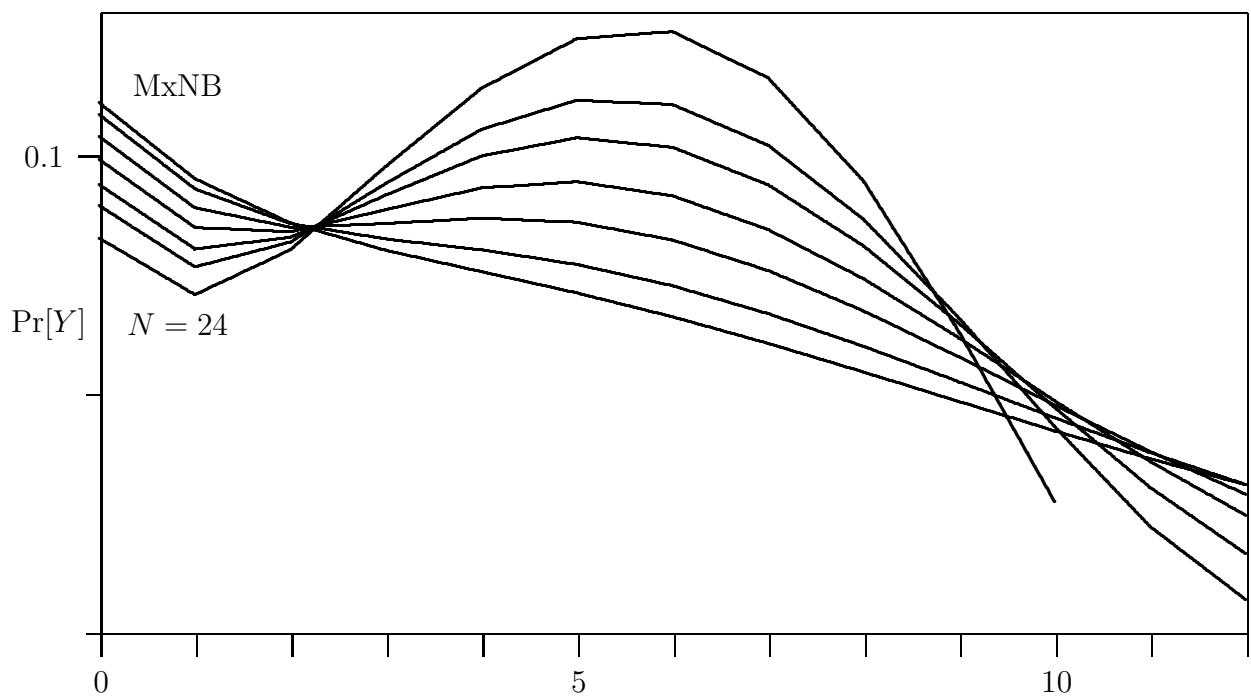


\caption { The maximum negative hypergeometric
distribution for $\,c=6\,$ and $\,m=N/3\,$ with $\,N=24,\; 27,\; 30, \; 36, \; 48, \; {\rm and \ }\; 96$.
Bimodal distributions generally appear whenever $\,c\,$ is moderately large and $\,N\,$ is small.
Examples of these conditions are given in Table~\ref{modetable}.
The limiting maximum negative binomial distribution (MxNB) has parameters $\,c=6\,$ and $\,p=1/3$ }

\bigskip
\beginpicture

\setcoordinatesystem units < .5in, 25in>
\setplotarea x from 0 to 12, y from 0  to  .13

\axis top   /
\axis right /

\axis bottom ticks out
    numbered from 0 to 10 by 5  
    unlabeled from 0 to 12 by 1  / 

\axis left ticks out
    numbered from .1 to .1 by .1  
    unlabeled short from 0 to .05 by .05 /

\put{$N=24$} at .8 .065
\put{MxNB}    at .8 .115
\put{$\Pr[Y]$} at -.55 .065

\setlinear

\plot   
0 .082918 1 .071072 2 .080495 3 .097994 4 .114326
5 .124720 6 .126204 7 .116496 8 .094861 9 .063241 10 .027667 /

\plot  
 0 .0897025 1 .076887 2 .082105 3 .094401  4 .105720 5 .111836
 6 .110904  7 .102373 8 .086834 9 .066159 10 .04341 11 .022474
12 .007179 /

\plot  
 0 .09410679  1 .08066297  2 .08310131 3 .09195402 4 .10014993
 5 .10394803  6 .10194903  7 .09410679 8 .08127405 9 .06501924
10 .04740986 11 .03067697 12 .01679929 /

\plot  
 0 .09936001  1 .08516572  2 .08424001 3 .08892001 4 .09344806
 5 .09474929  6 .09178099  7 .08471635 8 .07439377 9 .06199481
10 .04882091 11 .03610286 12 .02484000 /

\plot 
 0 .10416159  1 .08928137  2 .08520702 3 .08597465 4 .08712404
 5 .08624632  6 .08252798  7 .07612414 8 .06767634 9 .05800681
10 .04793603 11 .03817530 12 .02926773 /

\plot  
 0 .10876354  1 .09322589  2 .08595850 3 .08274266 4 .08046440
 5 .07740792  6 .07294153  7 .06709629 8 .06022936 9 .05280524
10 .04527182 11 .03799984 12 .03126094 /

\plot 
 0 .11124144  1 .09534980  2 .08608313 3 .08034907 4 .07589014
 5 .07143049  6 .06644389  7 .06086442 8 .05485730 9 .04866971
10 .04254666 11 .03669088 12 .03124914 /

\endpicture 

\end{figure}


\typeout{ Figure 3}
\begin{figure}[ht]
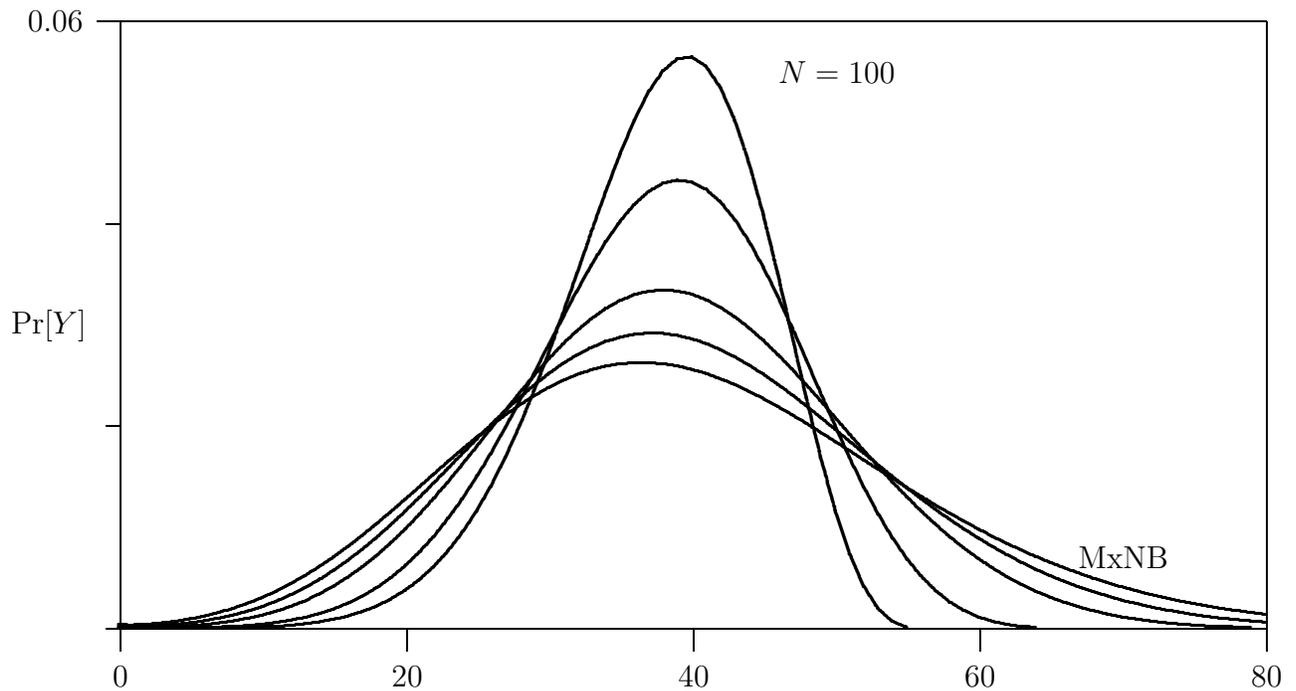


\caption{The maximum negative hypergeometric
distribution for $\,c=20\,$ and $\,m=N/4\,$ with $\,N=100, \; 120, 
\;200,\,$ and 400. A normal approximate distribution appears when the $\,(c,m,N)\,$ parameters are all large and $\,m/N\,$ is not close to zero, one, or 1/2.
This is proved in Lemma~4 }
\bigskip
\beginpicture

\setcoordinatesystem units < .075in, 53in>
\setplotarea x from 0 to 80, y from 0  to  .06

\axis top   /
\axis right /

\axis bottom ticks out
    numbered from 0 to 80 by 20  / 

\axis left ticks out
    numbered from  .06 to  .06 by  .06
    unlabeled short from 0 to .04 by .02 /

\put{$N=100$} at 50 .055
\put{MxNB}     at 70 .007
\put{$\Pr[Y]$} at -5 .03

\setlinear

\plot  
 0 .0003975  1 .0003786  2 .0004410  3 .0005637  4 .0007395  5 .0009682
 6 .0012533  7 .0016000  8 .0020137  9 .0024996 10 .0030619 11 .0037039
12 .0044273 13 .0052322 14 .0061171 15 .0070784 16 .0081106 17 .0092067
18 .0103575 19 .0115526 20 .0127801 21 .0140269 22 .0152793 23 .0165230
24 .0177434 25 .0189263 26 .0200578 27 .0211247 28 .0221149 29 .0230176 
30 .0238232 31 .0245239 32 .0251134 33 .0255872 34 .0259426 35 .0261785
36 .0262953 37 .0262953 38 .0261820 39 .0259601 40 .0256356 41 .0252153
42 .0247070 43 .0241187 44 .0234592 45 .0227374 46 .0219622 47 .0211428
48 .0202877 49 .0194056 50 .0185047 51 .0175925 52 .0166762 53 .0157624
54 .0148571 55 .0139657 56 .0130928 57 .0122427 58 .0114186 59 .0106237
60 .0098601 61 .0091297 62 .0084338 63 .0077733 64 .0071487 65 .0065600
66 .0060070 67 .0054891 68 .0050057 69 .0045557 70 .0041381 71 .0037516
72 .0033948 73 .0030662 74 .0027645 75 .0024881 76 .0022354 77 .0020049
78 .0017952 79 .0016048 80 .0014323 /

\plot  
 4 .0000126  5 .0000202  6 .0000317  7 .0000491  8 .0000747  9 .0001117
10 .0001646 11 .0002390 12 .0003421 13 .0004829 14 .0006727 15 .0009251
16 .0012564 17 .0016855 18 .0022343 19 .0029272 20 .0037911 21 .0048544
22 .0061466 23 .0076964 24 .0095307 25 .0116722 26 .0141372 27 .0169328
28 .0200542 29 .0234819 30 .0271784 31 .0310864 32 .0351268 33 .0391981
34 .0431771 35 .0469211 36 .0502726 37 .0530655 38 .0551341 39 .0563234
40 .0565022 41 .0555759 42 .0535001 43 .0502920 44 .0460393 45 .0409042
46 .0351197 47 .0289794 48 .0228163 49 .0169745 50 .0117719 51 .0074610
52 .0041911 53 .0019807 54 .0007112 55 .0001485  /

\plot  
 0 .00001337  1 .00001273  2 .00001847  3 .00002893  4 .00004502
 5 .00006880  6 .00010319  7 .00015205  8 .00022027  9 .00031395
10 .00044057 11 .00060908 12 .00083004 13 .00111561 14 .00147948
15 .00193677 16 .00250373 17 .00319734 18 .00403473 19 .00503257 
20 .00620615 21 .00756848 22 .00912921 23 .01089357 24 .01286118 
25 .01502513 26 .01737095 27 .01987598 28 .02250892 29 .02522978
30 .02799022 31 .03073436 32 .03340009 33 .03592085 34 .03822783 
35 .04025254 36 .04192973 37 .04320033 38 .04401445 39 .04433417 
40 .04413593 41 .04341239 42 .04217357 43 .04044717 44 .03827792
45 .03572605 46 .03286488 47 .02977748 48 .02655292 49 .0232819  
50 .02005248 51 .01694576 52 .01403219 53 .01136855 54 .00899583 
55 .00693832 56 .00520374 57 .00378454 58 .00266014 59 .00179996
60 .00116676 61 .00072022 62 .00042020 63 .00022951 64 .00011588  /

\plot  
 0 .00008  1 .00008  2 .00010  3 .00014  4 .00020  5 .00029  6 .00041 
 7 .00057  8 .00076  9 .00102 10 .00133 11 .00172 12 .00219 13 .00275
14 .00342 15 .00419 16 .00508 17 .00609 18 .00722 19 .00847 20 .00983
21 .01131 22 .01288 23 .01453 24 .01625 25 .01802 26 .01980 27 .02158
28 .02333 29 .02502 30 .02662 31 .02810 32 .02945 33 .03063 34 .03163
35 .03242 36 .03300 37 .03335 38 .03348 39 .03338 40 .03305 41 .03251
42 .03176 43 .03083 44 .02973 45 .02849 46 .02712 47 .02565 48 .02410
49 .02250 50 .02088 51 .01925 52 .01763 53 .01605 54 .01451 55 .01304
56 .01164 57 .01033 58 .00910 59 .00797 60 .00693 61 .00599 62 .00514
63 .00438 64 .00371 65 .00312 66 .00261 67 .00216 68 .00178 69 .00145
70 .00118 71 .00095 72 .00076 73 .00060 74 .00047 75 .00037 76 .00028
77 .00022 78 .00016 79 .00012  /

\plot  
 0  .00020  1  .00019  2  .00023  3  .00031  4  .00043  5  .00058
 6  .00078  7  .00104  8  .00135  9  .00172 10  .00217 11  .00271
12  .00333 13  .00404 14  .00484 15  .00575 16  .00674 17  .00783
18  .00901 19  .01027 20  .01159 21  .01297 22  .01440 23  .01585
24  .01730 25  .01875 26  .02017 27  .02154 28  .02284 29  .02406
30  .02518 31  .02618 32  .02705 33  .02779 34  .02838 35  .02881
36  .02909 37  .02921 38  .02917 39  .02898 40  .02865 41  .02818
42  .02758 43  .02687 44  .02605 45  .02514 46  .02416 47  .02311
48  .02201 49  .02087 50  .01971 51  .01854 52  .01736 53  .01620
54  .01505 55  .01393 56  .01285 57  .01180 58  .01080 59  .00985
60  .00895 61  .00810 62  .00731 63  .00657 64  .00588 65  .00525
66  .00467 67  .00414 68  .00366 69  .00323 70  .00283 71  .00248
72  .00216 73  .00188 74  .00163 75  .00141 76  .00121 77  .00104
78  .00089 79  .00076 80  .00065 /

\endpicture 
\end{figure}


\typeout{ Figure 4}
\begin{figure}[ht]
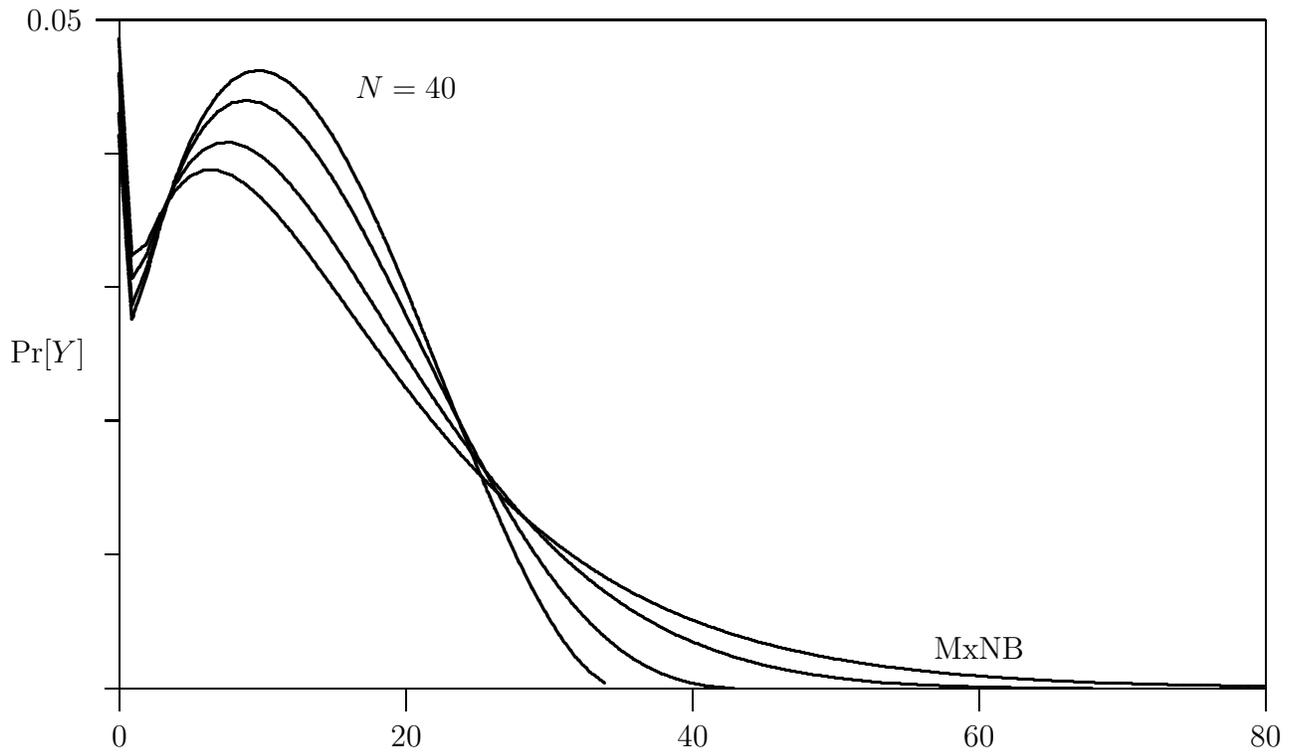

\caption{The maximum negative hypergeometric
distribution for $\,c=2\,$ and $\,m=N/10\,$ with $\,N=40,\; 50,\,$ and 100.
A gamma approximate distribution is proved in Lemma~2 when $\,m\,$ is much smaller than $\,N$. There is always a local mode at $\,Y=0$.  The asymptotic maximum negative binomial (MxNB) distribtion has parameters $\,c=2\,$ and $\,p=1/10$}
\bigskip
\beginpicture

\setcoordinatesystem units < .075in, 70in>
\setplotarea x from 0 to 80, y from 0  to  .05

\axis top   /
\axis right /
\axis bottom ticks out
    numbered from 0 to 80 by 20  / 
\axis left ticks out
    numbered from  .05 to  .05 by  .05
    unlabeled short from 0 to .04 by .01  / 

\put{$N=40$} at 20 .045
\put{MxNB}    at 60 .003
\put{$\Pr[Y]$} at -5 .025

\setlinear

\plot  
.1 .04136  1 .02757  2 .03074  3 .03466  4 .03799  5 .04070  6 .04283
 7 .04442  8 .04549  9 .04608 10 .04623 11 .04595 12 .04530 13 .04429
14 .04296 15 .04136 16 .03950 17 .03742 18 .03515 19 .03273 20 .03020
21 .02757 22 .02489 23 .02219 24 .01949 25 .01685 26 .01427 27 .01181
28 .00949 29 .00735 30 .00541 31 .00372 32 .00229 33 .00118 34 .00040 /

\plot   
.1 .04298  1 .02865  2 .03135  3 .03495  4 .03792  5 .04024  6 .04196 
 7 .04313  8 .04379  9 .04400 10 .04380 11 .04324 12 .04236 13 .04120
14 .03979 15 .03818 16 .03640 17 .03449 18 .03246 19 .03037 20 .02822
21 .02605 22 .02388 23 .02173 24 .01962 25 .01757 26 .01560 27 .01372
28 .01193 29 .01027 30 .00872 31 .00730 32 .00601 33 .00485 34 .00384
35 .00295 36 .00220 37 .00158 38 .00108 39 .00069 40 .00040 41 .00020
42 .00008 43 .00002 /

\plot  
.1 .04596  1 .03064  2 .03233  3 .03531  4 .03764  5 .03927  6 .04030
 7 .04079  8 .04084  9 .04050 10 .03984 11 .03892 12 .03777 13 .03645
14 .03500 15 .03344 16 .03181 17 .03014 18 .02844 19 .02674 20 .02505
21 .02339 22 .02176 23 .02019 24 .01866 25 .01720 26 .01581 27 .01449
28 .01323 29 .01205 30 .01094 31 .00991 32 .00894 33 .00805 34 .00722
35 .00646 36 .00576 37 .00512 38 .00453 39 .00400 40 .00352 41 .00309
42 .00270 43 .00235 44 .00204 45 .00176 46 .00151 47 .00130 48 .00111
49 .00094 50 .00079 51 .00067 52 .00056 53 .00046 54 .00038 55 .00032
56 .00026 57 .00021 58 .00017 59 .00013 60 .00011 61 .00008 62 .00006
63 .00005 64 .00004 65 .00003 66 .00002 67 .00002 68  .00001 /

\plot 
.1   .0486        1   .0324        2   .03321       3   .035478
 4   .03720654    5   .0382644     6   .03874212    7   .03874206
 8   .03835463    9   .03765727   10   .03671584   11   .03558612
12   .03431519   13   .03294258   14   .03150134   15   .03001893
16   .02851798   17   .02701703   18   .0255311    19   .02407218
20   .02264973   21   .02127105   22   .01994161   23   .01866535
24   .01744492   25   .01628193   26   .01517708   27   .01413039
28   .01314126   29   .01220865   30   .01133116   31   .01050707
32   .009734494  33   .00901136   34   .008335508  35   .007704713
36   .007116721  37   .006569281  38   .006060162  39   .005587174
40   .005148182  41   .004741116  42   .004363982  43   .004014863
44   .003691929  45   .003393432  46   .003117716  47   .002863208
48   .002628425  49   .002411967  50   .002212516  51   .002028835
52   .001859766  53   .001704222  54   .001561189  55   .00142972
56   .001308933  57   .001198007  58   .001096176  59   .001002732
60   .0009170144 61   .0008384131 62   .000766362  63   .000700337
64   .0006398533 65   .000584463  66   .0005337523 67   .000487339
68   .0004448709 69   .000406023  70   .000370496  71   .0003380142
72   .0003083237 73   .0002811912 74   .000256402  75   .0002337587
76   .0002130801 77   .0001941996 78   .0001769643 79   .0001612342
80   .0001468804    /

\endpicture 

\end{figure}


\typeout{ Figure 5}
\begin{figure}[ht]
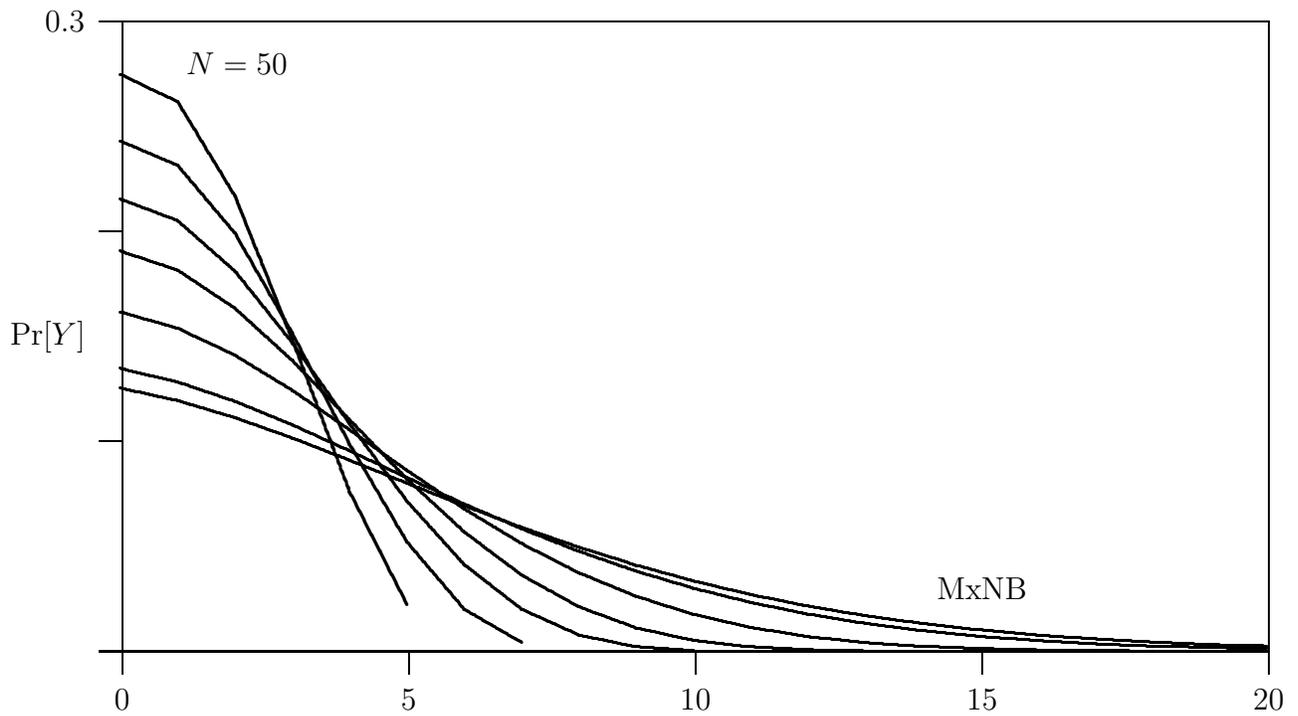

\caption{The maximum negative hypergeometric
distribution for $\,c=20\,$ and $\,m=N/2\,$ for
 $\,N=50,\; 54,\; 60,\; 70,\; 100,\,$ and 300. A half normal approximate distribution appears when $\,c\,$ and $\,N\,$ are large and $\,m=N/2$.  The limiting distribution is proved in Lemma~3.}
\bigskip
\beginpicture

\setcoordinatesystem units < .3in, 11in>
\setplotarea x from 0 to 20, y from 0  to  .3

\axis top   /
\axis right /

\axis bottom ticks out
    numbered from 0 to 20 by 5  / 

\axis left ticks out
    numbered from  .3 to  .3 by  .3
    unlabeled from 0 to .2 by .1   / 

\put{$N=50$} at 2 .28
\put{MxNB}    at 15 .03
\put{$\Pr[Y]$} at -1.3 .15

\setlinear
\plot    
0  .27479 1  .26171 2  .21677 3  .14845 4  .07599 5  .02229  /

\plot   
0 .24299 1 .23142 2 .19905 3 .15145 4 .09867 5 .05210 6 .02003 7 .00426 /

\plot  
0 .21534 1 .20509 2 .18105  3 .14693 4 .10840 5 .07154 6 .04127 
7 .02009 8 .00778 9 .00214 10 .00031  /

\plot  
 0  .19059  1 .18152  2 .16331  3 .13846 4  .11025  5 .08209  6 .05683
 7  .03631  8 .02120  9 .01116 10 .00521 11 .00210 12 .00070 13 .00018
14 .00003 /

\plot  
 0  .16158335608  1  .15388891055  2  .14096603594  3  .12426990874
 4  .10546590939  5  .08618071453  6  .06779951318  7  .05133790298
 8  .03739641046  9  .02618740680 10  .01761231477 11  .01136278372
12  .00702202642 13  .00414937925 14  .00233954362 15  .00125550664
16  .00063937838 17  .00030790703 18  .00013963226 18  .00005933092
20  .00002347913 /

\plot  
 0  .134652457 1  .128240435 2  .119035392 3  .107841973 4  .095480709
 5  .082710164 6  .070172651 7  .058364654 8  .047629315 9  .038166019
10  .030051305 11  .023265527 12  .017720693 13  .013286189
14  .009810390 15  .007137315 16  .005118313 17  .003619330 18  .002524594
19  .001737631 20  .001180472 /

\plot 
 0  .125370688  1  .119400655  2  .111259701  3  .101584945  4  .091003179
 5  .080082798  6  .069302421  7  .059035396  8  .049547564  9  .041004881
10  .033487319 11  .027005903 12  .021520329 13  .016955411 14  .013215246 15  .010194619 16  .007787556 17  .005893286 18  .004419964 19  .003286640 20  .002423897 /

\endpicture

\end{figure}


\typeout{ Figure 6}
\begin{figure}[ht]
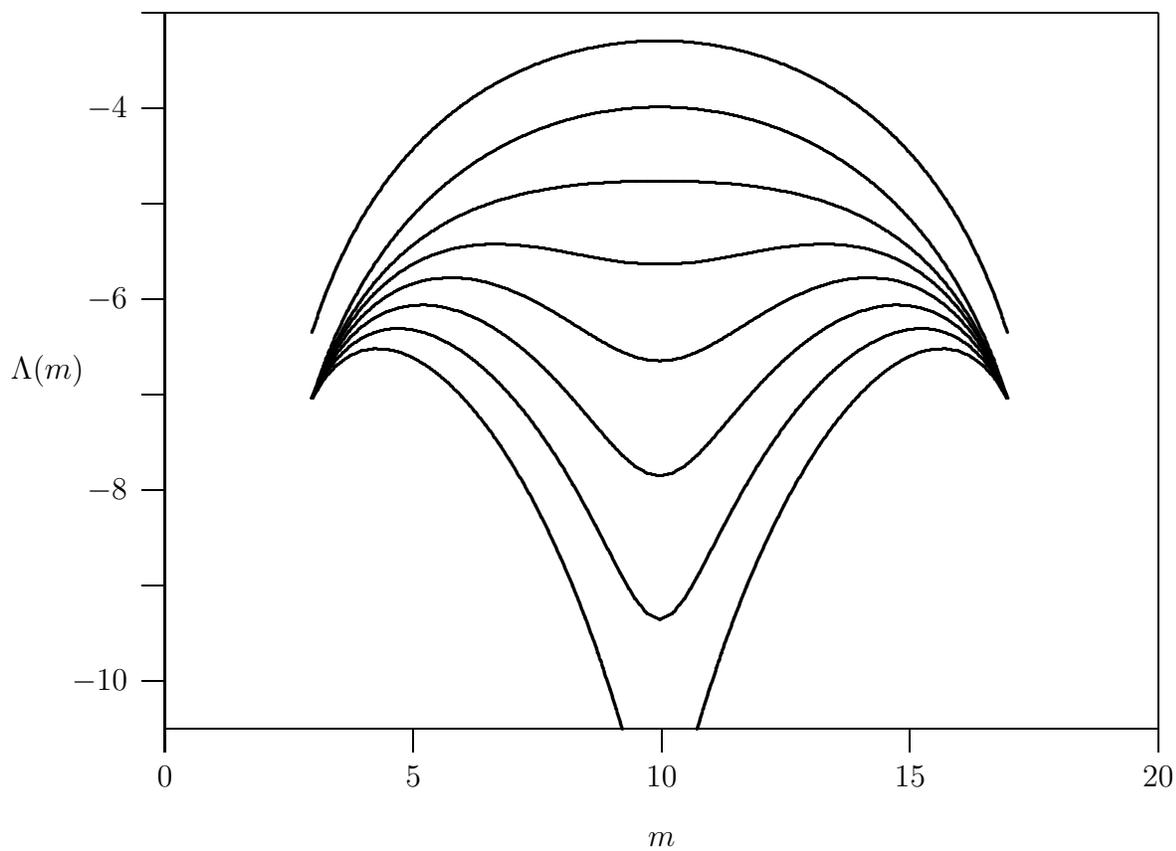

\caption{The maximum negative hypergeometric log-likelihood kernel function  $\,\Lambda(m)\,$ for parameter values $\,N=20\,$ and $\,c=3\,$ and observed values $\,y=0,\ldots,7$, from top to bottom, respectively }
\bigskip
\beginpicture

\setcoordinatesystem units < .26in, .5in>
\setplotarea x from 0 to 20, y from -10.5  to  -3

\axis top   /
\axis right /

\axis bottom 
     label {$m$}
     ticks out
     numbered from 0 to 20 by 5  / 

\axis left 
     label {$\Lambda(m)\!\!\!\!\!\!\!\!\!\!\!$}
     ticks out
     numbered from  -10 to -3 by 2
     unlabeled from -10 to -3 by 1   / 

\setlinear

\plot  
  3.00 -6.345636  3.25 -5.972038  3.50 -5.658380  3.75 -5.388846
  4.00 -5.153498  4.25 -4.945622  4.50 -4.760439  4.75 -4.594413
  5.00 -4.444847  5.25 -4.309633  5.50 -4.187093  5.75 -4.075867
  6.00 -3.974843  6.25 -3.883096  6.50 -3.799855  6.75 -3.724470
  7.00 -3.656389  7.25 -3.595147  7.50 -3.540344  7.75 -3.491641
  8.00 -3.448750  8.25 -3.411426  8.50 -3.379463  8.75 -3.352690
  9.00 -3.330967  9.25 -3.314183  9.50 -3.302254  9.75 -3.295120
 10.00 -3.292746 10.25 -3.295120 10.50 -3.302254 10.75 -3.314183
 11.00 -3.330967 11.25 -3.352690 11.50 -3.379463 11.75 -3.411426
 12.00 -3.448750 12.25 -3.491641 12.50 -3.540344 12.75 -3.595147
 13.00 -3.656389 13.25 -3.724470 13.50 -3.799855 13.75 -3.883096
 14.00 -3.974843 14.25 -4.075867 14.50 -4.187093 14.75 -4.309633
 15.00 -4.444847 15.25 -4.594413 15.50 -4.760439 15.75 -4.945622
 16.00 -5.153498 16.25 -5.388846 16.50 -5.658380 16.75 -5.972038
 17.00 -6.345636   /

\plot 
  3.00 -7.038784  3.25 -6.665185  3.50 -6.351527  3.75 -6.081993
  4.00 -5.846645  4.25 -5.638769  4.50 -5.453586  4.75 -5.287560
  5.00 -5.137994  5.25 -5.002780  5.50 -4.880240  5.75 -4.769015
  6.00 -4.667990  6.25 -4.576243  6.50 -4.493002  6.75 -4.417617
  7.00 -4.349536  7.25 -4.288294  7.50 -4.233491  7.75 -4.184788
  8.00 -4.141897  8.25 -4.104573  8.50 -4.072610  8.75 -4.045837
  9.00 -4.024114  9.25 -4.007330  9.50 -3.995401  9.75 -3.988267
 10.00 -3.985893 10.25 -3.988267 10.50 -3.995401 10.75 -4.007330
 11.00 -4.024114 11.25 -4.045837 11.50 -4.072610 11.75 -4.104573
 12.00 -4.141897 12.25 -4.184788 12.50 -4.233491 12.75 -4.288294
 13.00 -4.349536 13.25 -4.417617 13.50 -4.493002 13.75 -4.576243
 14.00 -4.667990 14.25 -4.769015 14.50 -4.880240 14.75 -5.002780
 15.00 -5.137994 15.25 -5.287560 15.50 -5.453586 15.75 -5.638769
 16.00 -5.846645 16.25 -6.081993 16.50 -6.351527 16.75 -6.665185
 17.00 -7.038784  /

\plot 
  3.00 -7.038784  3.25 -6.703691  3.50 -6.428598  3.75 -6.197632
  4.00 -6.000796  4.25 -5.831307  4.50 -5.684313  4.75 -5.556194
  5.00 -5.444161  5.25 -5.346006  5.50 -5.259941  5.75 -5.184488
  6.00 -5.118407  6.25 -5.060636  6.50 -5.010259  6.75 -4.966473
  7.00 -4.928570  7.25 -4.895919  7.50 -4.867955  7.75 -4.844171
  8.00 -4.824115  8.25 -4.807382  8.50 -4.793614  8.75 -4.782500
  9.00 -4.773773  9.25 -4.767216  9.50 -4.762656  9.75 -4.759970
 10.00 -4.759083 10.25 -4.759970 10.50 -4.762656 10.75 -4.767216
 11.00 -4.773773 11.25 -4.782500 11.50 -4.793614 11.75 -4.807382
 12.00 -4.824115 12.25 -4.844171 12.50 -4.867955 12.75 -4.895919
 13.00 -4.928570 13.25 -4.966473 13.50 -5.010259 13.75 -5.060636
 14.00 -5.118407 14.25 -5.184488 14.50 -5.259941 14.75 -5.346006
 15.00 -5.444161 15.25 -5.556194 15.50 -5.684313 15.75 -5.831307
 16.00 -6.000796 16.25 -6.197632 16.50 -6.428598 16.75 -6.703691
 17.00 -7.038784  /

\plot 
  3.00 -7.038784  3.25 -6.723516  3.50 -6.469482  3.75 -6.260887
  4.00 -6.087807  4.25 -5.943537  4.50 -5.823297  4.75 -5.723534
  5.00 -5.641520  5.25 -5.575096  5.50 -5.522506  5.75 -5.482280
  6.00 -5.453156  6.25 -5.434013  6.50 -5.423821  6.75 -5.421606
  7.00 -5.426409  7.25 -5.437260  7.50 -5.453154  7.75 -5.473029
  8.00 -5.495756  8.25 -5.520138  8.50 -5.544917  8.75 -5.568807
  9.00 -5.590535  9.25 -5.608901  9.50 -5.622857  9.75 -5.631583
 10.00 -5.634551 10.25 -5.631583 10.50 -5.622857 10.75 -5.608901
 11.00 -5.590535 11.25 -5.568807 11.50 -5.544917 11.75 -5.520138
 12.00 -5.495756 12.25 -5.473029 12.50 -5.453154 12.75 -5.437260
 13.00 -5.426409 13.25 -5.421606 13.50 -5.423821 13.75 -5.434013
 14.00 -5.453156 14.25 -5.482280 14.50 -5.522506 14.75 -5.575096
 15.00 -5.641520 15.25 -5.723534 15.50 -5.823297 15.75 -5.943537
 16.00 -6.087807 16.25 -6.260887 16.50 -6.469482 16.75 -6.723516
 17.00 -7.038784 /

\plot 
  3.00 -7.038784  3.25 -6.746705  3.50 -6.516241  3.75 -6.331661
  4.00 -6.183117  4.25 -6.063994  4.50 -5.969612  4.75 -5.896542
  5.00 -5.842191  5.25 -5.804558  5.50 -5.782064  5.75 -5.773443
  6.00 -5.777652  6.25 -5.793811  6.50 -5.821143  6.75 -5.858925
  7.00 -5.906433  7.25 -5.962879  7.50 -6.027344  7.75 -6.098689
  8.00 -6.175455  8.25 -6.255747  8.50 -6.337124  8.75 -6.416515
  9.00 -6.490218  9.25 -6.554044  9.50 -6.603687  9.75 -6.635291
 10.00 -6.646152 10.25 -6.635291 10.50 -6.603687 10.75 -6.554044
 11.00 -6.490218 11.25 -6.416515 11.50 -6.337124 11.75 -6.255747
 12.00 -6.175455 12.25 -6.098689 12.50 -6.027344 12.75 -5.962879
 13.00 -5.906433 13.25 -5.858925 13.50 -5.821143 13.75 -5.793811
 14.00 -5.777652 14.25 -5.773443 14.50 -5.782064 14.75 -5.804558
 15.00 -5.842191 15.25 -5.896542 15.50 -5.969612 15.75 -6.063994
 16.00 -6.183117 16.25 -6.331661 16.50 -6.516241 16.75 -6.746705
 17.00 -7.038784  /

\plot 
  3.00 -7.038784  3.25 -6.771966  3.50 -6.567472  3.75 -6.409584
  4.00 -6.288478  4.25 -6.197559  4.50 -6.132182  4.75 -6.088954
  5.00 -6.065334  5.25 -6.059390  5.50 -6.069631  5.75 -6.094907
  6.00 -6.134327  6.25 -6.187211  6.50 -6.253037  6.75 -6.331408
  7.00 -6.422009  7.25 -6.524553  7.50 -6.638710  7.75 -6.763999
  8.00 -6.899612  8.25 -7.044151  8.50 -7.195224  8.75 -7.348874
  9.00 -7.498882  9.25 -7.636162  9.50 -7.748819  9.75 -7.823736
 10.00 -7.850125 10.25 -7.823736 10.50 -7.748819 10.75 -7.636162
 11.00 -7.498882 11.25 -7.348874 11.50 -7.195224 11.75 -7.044151
 12.00 -6.899612 12.25 -6.763999 12.50 -6.638710 12.75 -6.524553
 13.00 -6.422009 13.25 -6.331408 13.50 -6.253037 13.75 -6.187211
 14.00 -6.134327 14.25 -6.094907 14.50 -6.069631 14.75 -6.059390
 15.00 -6.065334 15.25 -6.088954 15.50 -6.132182 15.75 -6.197559
 16.00 -6.288478 16.25 -6.409584 16.50 -6.567472 16.75 -6.771966
 17.00 -7.038784  /

\plot 
  3.00 -7.038784  3.25 -6.800162  3.50 -6.624656  3.75 -6.496611
  4.00 -6.406261  4.25 -6.347079  4.50 -6.314487  4.75 -6.305164
  5.00 -6.316649  5.25 -6.347087  5.50 -6.395080  5.75 -6.459574
  6.00 -6.539792  6.25 -6.635189  6.50 -6.745412  6.75 -6.870287
  7.00 -7.009796  7.25 -7.164071  7.50 -7.333372  7.75 -7.518048
  8.00 -7.718447  8.25 -7.934715  8.50 -8.166354  8.75 -8.411312
  9.00 -8.664185  9.25 -8.913040  9.50 -9.135113  9.75 -9.294990
 10.00 -9.354203 10.25 -9.294990 10.50 -9.135113 10.75 -8.913040
 11.00 -8.664185 11.25 -8.411312 11.50 -8.166354 11.75 -7.934715
 12.00 -7.718447 12.25 -7.518048 12.50 -7.333372 12.75 -7.164071
 13.00 -7.009796 13.25 -6.870287 13.50 -6.745412 13.75 -6.635189
 14.00 -6.539792 14.25 -6.459574 14.50 -6.395080 14.75 -6.347087
 15.00 -6.316649 15.25 -6.305164 15.50 -6.314487 15.75 -6.347079
 16.00 -6.406261 16.25 -6.496611 16.50 -6.624656 16.75 -6.800162
 17.00 -7.038784  /

\plot  
  3.00  -7.038784  3.25  -6.831895  3.50  -6.689179  3.75  -6.595042
  4.00  -6.539792  4.25  -6.516985  4.50  -6.522135  4.75  -6.552031
  5.00  -6.604331  5.25  -6.677322  5.50  -6.769761  5.75  -6.880777
  6.00  -7.009796  6.25  -7.156503  6.50  -7.320810  6.75  -7.502844
  7.00  -7.702943  7.25  -7.921669  7.50  -8.159834  7.75  -8.418546
  8.00  -8.699277  8.25  -9.003939  8.50  -9.334924  8.75  -9.694885
  9.00 -10.085571  9.25 -10.503555 /


\plot
 10.75 -10.503555
 11.00 -10.085571 11.25  -9.694885 11.50  -9.334924 11.75  -9.003939
 12.00  -8.699277 12.25  -8.418546 12.50  -8.159834 12.75  -7.921669
 13.00  -7.702943 13.25  -7.502844 13.50  -7.320810 13.75  -7.156503
 14.00  -7.009796 14.25  -6.880777 14.50  -6.769761 14.75  -6.677322
 15.00  -6.604331 15.25  -6.552031 15.50  -6.522135 15.75  -6.516985
 16.00  -6.539792 16.25  -6.595042 16.50  -6.689179 16.75  -6.831895
 17.00  -7.038784 /

\endpicture
\end{figure}
\typeout{end of figures}

\end{document}